\documentclass[reqno,12pt]{amsart}
\usepackage{amsmath,amsrefs,amssymb}

\numberwithin{equation}{section}

\DeclareMathOperator{\divergence}{div}
\DeclareMathOperator{\loc}{loc}
\DeclareMathOperator{\Tr}{Tr}
\DeclareMathOperator{\bigO}{O}
\DeclareMathOperator{\smallo}{o}
\newcommand{\R}{\mathbb{R}}

\newcommand{\<}{\left<}
\renewcommand{\>}{\right>}
\renewcommand{\[}{\left[}
\renewcommand{\]}{\right]}
\renewcommand{\(}{\left(}
\renewcommand{\)}{\right)}

\newtheorem{theorem}{Theorem}[section]
\newtheorem{lemma}[theorem]{Lemma}

\begin{document}

\title[Positive solutions to the critical $p$-Laplace equation in $\R^n$]{A note on the classification of positive solutions to the critical $p$-Laplace equation in $\R^n$}

\author{J\'er\^ome V\'etois}

\address{J\'er\^ome V\'etois, Department of Mathematics and Statistics, McGill University, 805 Sherbrooke Street West, Montreal, Quebec H3A 0B9, Canada}
\email{jerome.vetois@mcgill.ca}

\thanks{To appear in {\it Advanced Nonlinear Studies}.}

\thanks{The author was supported by the NSERC Discovery Grant RGPIN-2022-04213.}

\date{February 19, 2024}

\begin{abstract}
In this note, we obtain a classification result for positive solutions to the critical $p$-Laplace equation in $\R^n$ with $n\ge4$ and $p>p_n$ for some number $p_n\in\(\frac{n}{3},\frac{n+1}{3}\)$ such that $p_n\sim\frac{n}{3}+\frac{1}{n}$, which improves upon a similar result obtained by Ou~\cite{Ou} under the condition $p\ge\frac{n+1}{3}$.
\end{abstract}

\maketitle

\section{Introduction and main result}

We consider positive, weak solutions $u\in W^{1,p}_{\loc}\(\R^n\)\cap L^\infty_{\loc}\(\R^n\)$ to the critical $p$-Laplace equation 
\begin{equation}\label{IntroEq1}
-\Delta_p u=u^{p^*-1}\quad\text{in }\R^n,
\end{equation}
where $n\ge2$, $1<p<n$, $\Delta_p:=\divergence\(\left|\nabla u\right|^{p-2}\nabla u\)$ is the $p$-Laplace operator and $p^*:=np/\(n-p\)$ is the critical Sobolev exponent.

\smallskip
Well-known solutions to \eqref{IntroEq1} are the functions
\begin{equation}\label{IntroEq2}
u_{\mu,x_0}\(x\):=\(\frac{n^{\frac{1}{p}}\(\frac{n-p}{p-1}\)^{\frac{p-1}{p}}\mu^{\frac{1}{p-1}}}{\mu^{\frac{p}{p-1}}+\left|x-x_0\right|^{\frac{p}{p-1}}}\)^{\frac{n-p}{p}}\quad\forall x\in\R^n,
\end{equation}
where $\mu>0$ and $x_0\in\R^n$. As was shown by Rodemich~\cite{Rod}, Aubin~\cite{Aub} and Talenti~\cite{Tal}, these functions realize the equality in the optimal Sobolev inequality in $\R^n$. Guedda and V\'eron~\cite{GueVer} obtained that the functions defined in \eqref{IntroEq2} are the only positive, radially symmetric solutions to \eqref{IntroEq1}. In the case where $p=2$, Caffarelli, Gidas and Spruck~\cite{CafGidSpr} (see also Chen and Li~\cite{ChenLi}) used the moving plane method to obtain that these functions are in fact the only positive solutions of \eqref{IntroEq1}. This classification result was later extended by Damascelli and Ramaswamy~\cite{DamRam} to the case of solutions with sufficiently fast decay at infinity with $1<p<2$, and in a series of papers by Damascelli, Merch\'an, Montoro and Sciunzi~\cite{DamMerMonSci}, V\'etois~\cite{Vet} and Sciunzi~\cite{Sci} to the case of solutions in $D^{1,p}\(\R^n\)$ for all $p\in\(1,n\)$. We mention in passing that a similar classification result was also obtained by Esposito~\cite{Esp} for solutions with finite mass of the critical $n$-Laplace equation, in which case the nonlinearity is of exponential type.

\smallskip
More recently, Ciraolo, Figalli and Roncoroni~\cite{CirFigRon} used a strategy based on integral estimates to extend the classification of positive $D^{1,p}$-solutions to a class of anisotropic $p$-Laplace-type equations in convex cones (see also the survey article by Roncoroni~\cite{Ron} on this topic). In the case where $p=2$, this type of approach can be traced back to the work of Obata~\cite{Oba} on the conformal transformations of the sphere. An approach of this type was then used by Catino, Monticelli and Roncoroni~\cite{CatMonRon} to obtain new classification results for positive, weak solutions to \eqref{IntroEq1} which are not a priori in $D^{1,p}\(\R^n\)$. In particular, Catino, Monticelli and Roncoroni~\cite{CatMonRon} managed to obtain the complete classification of positive, weak solutions to \eqref{IntroEq1} in the case where $n=2$ or [$n=3$ and $3/2<p<2$]. The method was recently improved by Ou~\cite{Ou} who managed to extend this result to the case where $n\ge3$ and $p\ge\(n+1\)/3$.

\smallskip
In this note, we obtain the following extension of Catino, Monticelli and Roncoroni~\cite{CatMonRon} and Ou's~\cite{Ou} results:

\begin{theorem}\label{Th}
Assume that $n\ge4$ and $p_n<p<n$, where
$$p_n:=\left\{\begin{aligned}
&\frac{8}{5}&&\text{if }n=4\\
&\frac{4n+3-\sqrt{4n^2+12n-15}}{6}&&\text{if }n\ge5.
\end{aligned}\right.$$
Then every positive, weak solution $u\in W^{1,p}_{\loc}\(\R^n\)\cap L^\infty_{\loc}\(\R^n\)$ to \eqref{IntroEq1} is of the form \eqref{IntroEq2}, i.e. $u\equiv u_{\mu,x_0}$ for some $\mu>0$ and $x_0\in\R^n$.
\end{theorem}

It is easy to see that
$$\frac{n}{3}<p_n<\frac{n+1}{3}\quad\forall n\ge4$$
and
$$p_n\sim\frac{n}{3}+\frac{1}{n}\quad\text{as }n\to\infty.$$ 

\smallskip
The main difficulty in our proof in the case where $p<\(n+1\)/3$ is to obtain a priori integral estimates with an exponent on the gradient which is larger than $p$. This can be seen for example by looking at the formula \eqref{ProofEq13} in our proof, where the exponent on the function g (defined in \eqref{ProofEq2}) is less than $1$ for small $\varepsilon>0$ if and only if $p>\(n+1\)/3$. While the former case can be achieved by using some rather straightforward estimates (see Lemma~\ref{Lem1}), the case where $p_n<p<\(n+1\)/3$ requires a little more work. In this case, by using the integral identity in Lemma~\ref{Lem3}, we manage to obtain the key estimate \eqref{ProofEq26}, which compares two integrals with different exponents on the gradient and from which we manage to derive our classification result. The case where $p\le p_n$ remains open. In this case, the exponent on the gradient in the right-hand side of \eqref{ProofEq26} becomes too large for us to conclude. The situation appears to be even more problematic when $p<n/3$ since the exponent on the gradient in the right-hand side of \eqref{ProofEq26} then becomes greater than the exponent in the left-hand side.

\section{Proof of Theorem~\ref{Th}}

Let $u\in W^{1,p}_{\loc}\(\R^n\)\cap L^\infty_{\loc}\(\R^n\)$ be a positive, weak solution of \eqref{IntroEq1}. Results by DiBenedetto~\cite{DiB} and Tolksdorf~\cite{Tol} give that $u\in C^{1,\alpha}_{\loc}\(\R^n\)$ for some $\alpha\in\(0,1\)$. Furthermore, as was shown by Antonini, Ciraolo and Farina~\cite{AntCirFar} (see also the references therein for previous results), the critical set $Z:=\left\{x\in\R^n:\,\left|\nabla u\(x\)\right|=0\right\}$ has measure zero, $u\in W^{2,2}_{\loc}\(\R^n\backslash Z\)$, $\left|\nabla u\right|^{p-2}\nabla u\in W^{1,2}_{\loc}\(\R^n\)$ and $\left|\nabla u\right|^{p-2}\nabla^2 u\in L^2_{\loc}\(\R^n\)$.

\smallskip
Following the approach developped by Catino, Monticelli and Roncoroni~\cite{CatMonRon} and Ou~\cite{Ou} (see also the previous work by Ciraolo, Figalli and Roncoroni~\cite{CirFigRon}), we define the function 
\begin{equation}\label{ProofEq1}
v:=u^{-\frac{p}{n-p}}.
\end{equation}
The equation \eqref{IntroEq1} can then be rewritten as
\begin{equation}\label{ProofEq2}
\Delta_pv=g:=\frac{n\(p-1\)}{p}v^{-1}\left|\nabla v\right|^p+\(\frac{p}{n-p}\)^{p-1}v^{-1}\quad\text{in }\R^n.
\end{equation}
Furthermore, it follows from the above-mentioned regularity properties of $u$ that $v\in C^{1,\alpha}_{\loc}\(\R^n\)\cap W^{2,2}_{\loc}\(\R^n\backslash Z\)$, $\left|\nabla v\right|^{p-2}\nabla v\in W^{1,2}_{\loc}\(\R^n\)$ and $\left|\nabla v\right|^{p-2}\nabla^2 v\in L^2_{\loc}\(\R^n\)$.

\smallskip
We now state some preliminary results, starting with the following lemma, of which more or less general versions can be found in either of the work by Serrin and Zou~\cite{SerZhou}*{Lemma~2.4}, Catino, Monticelli and Roncoroni~\cite{CatMonRon}*{Lemma~5.1} and Ou~\cite{Ou}*{Lemma~3.1}:

\begin{lemma}\label{Lem1}
Let $n\ge2$, $p\in\(1,n\)$, $r\in\[0,p\]$, $q<\(np-n+p\)/p$, $R>1$, $u\in W^{1,p}_{\loc}\(\R^n\)\cap L^\infty_{\loc}\(\R^n\)$ be a positive, weak solution of \eqref{IntroEq1} and $v$ be the function defined in \eqref{ProofEq1}. Then 
\begin{equation}\label{Lem1Eq1}
\int_{B_R\(0\)}v^{-q}\left|\nabla v\right|^r\le C\left\{\begin{aligned}&R^{n-q}&&\text{if }r\le q<\frac{np-n+p}{p}\\&R^{n-\frac{pq-r}{p-1}}&&\text{if }q<r\end{aligned}\right.
\end{equation}
for some constant $C=C\(n,p,q,r\)>0$.
\end{lemma}

\proof[Proof of Lemma~\ref{Lem1}]
Let $r\in\[0,p\]$ and $q<\(np-n+p\)/p$. We refer to Ou~\cite{Ou}*{Lemma~3.1} for the proof of \eqref{Lem1Eq1} when [$q\ge0$ and $r=0$] or [$q\ge p$ and $r=p$]. In the case where $q\ge r$ and $0<r<p$, H\"older's inequality gives
\begin{equation}\label{Lem1Eq2}
\int_{\R^n}v^{-q}\left|\nabla v\right|^r\le\(\int_{\R^n}v^{-q-\sigma\(p-r\)}\left|\nabla v\right|^p\)^{\frac{r}{p}}\(\int_{\R^n}v^{-q+\sigma r}\)^{\frac{p-r}{p}},
\end{equation}
where 
$$\sigma:=\max\(\frac{p-q}{p-r},0\),$$
so that
\begin{equation}\label{Lem1Eq3}
q+\sigma\(p-r\)=\max\(p,q\)\in\[p,\frac{np-n+p}{p}\)
\end{equation}
and
\begin{equation}\label{Lem1Eq4}
q-\sigma r=\min\(\frac{p\(q-r\)}{p-r},q\)\in\[0,\frac{np-n+p}{p}\).
\end{equation}
It follows from \eqref{Lem1Eq3} and \eqref{Lem1Eq4} (together with the above-mentioned proof by Ou~\cite{Ou}*{Lemma~3.1}) that
\begin{equation}\label{Lem1Eq5}
\int_{\R^n}v^{-q-\sigma\(p-r\)}\left|\nabla v\right|^p\le CR^{n-q-\sigma\(p-r\)}
\end{equation}
and
\begin{equation}\label{Lem1Eq6}
\int_{\R^n}v^{-q+\sigma r}\le CR^{n-q+\sigma r}
\end{equation}
for some constant $C=C\(n,p,q,r\)>0$. By combining \eqref{Lem1Eq2}, \eqref{Lem1Eq5} and \eqref{Lem1Eq6}, we then obtain
$$\int_{\R^n}v^{-q}\left|\nabla v\right|^r\le C\(R^{n-q-\sigma\(p-r\)}\)^{\frac{r}{p}}\(R^{n-q+\sigma r}\)^{\frac{p-r}{p}}=CR^{n-q}$$
for some constant $C=C\(n,p,q,r\)>0$. We now consider the case where $q<r$ and $0\le r\le p$. In this case, by observing that $\Delta_pu\le0$ in $\R^n$, we obtain (see Serrin and Zou~\cite{SerZhou}*{Lemma~2.3})
\begin{equation}\label{Lem1Eq7}
u\(x\)\ge C\left|x\right|^{-\frac{n-p}{p-1}},\text{ i.e. }v\(x\)\le C^{-\frac{p}{n-p}}\left|x\right|^{\frac{p}{p-1}}\quad\forall x\in\R^n\backslash B_1\(0\)
\end{equation}
for some constant $C=C\(n,p\)>0$. It follows from \eqref{Lem1Eq7} that
\begin{align*}
\int_{B_R\(0\)}v^{-q}\left|\nabla v\right|^r&\le C^{-\frac{p\(r-q\)}{n-p}}R^{\frac{p\(r-q\)}{p-1}}\int_{B_R\(0\)}v^{-r}\left|\nabla v\right|^r\\
&\le C^\prime R^{\frac{p\(r-q\)}{p-1}+n-r}\\
&=C^\prime R^{n-\frac{pq-r}{p-1}},
\end{align*}
for some constant $C^\prime=C^\prime\(n,p,q,r\)>0$. This ends the proof of Lemma~\ref{Lem1}. 
\endproof

Next, we state the following lemma obtained by Ou~\cite{Ou}*{Proposition~2.3}, which extends a previous result by Catino, Monticelli and Roncoroni~\cite{CatMonRon}*{Proposition~2.2} (see also Serrin and Zou~\cite{SerZhou}*{Proposition~6.2}):

\begin{lemma}
Let $n\ge2$, $p\in\(1,n\)$, $m\in\R$, $u\in W^{1,p}_{\loc}\(\R^n\)\cap L^\infty_{\loc}\(\R^n\)$ be a positive, weak solution of \eqref{IntroEq1}, $v$ and $g$ be the functions defined in \eqref{ProofEq1} and \eqref{ProofEq2}, and $\varphi$ be a smooth, nonnegative function with compact support in $\R^n$. Then 
\begin{multline}\label{Lem2Eq1}
\int_{\R^n}\varphi v^{1-n}g^m\Tr\(E^2\)+nm\int_{\R^n}\varphi v^{-n}g^{m-1}\left|\nabla v\right|^{p-2}\<E^2\nabla v,\nabla v\>\\
\le-\int_{\R^n}v^{1-n}g^m\left|\nabla v\right|^{p-2}\<E\nabla v,\nabla\varphi\>,
\end{multline}
where $E=\(E_{ij}\)_{1\le i,j\le n}$ is the matrix-valued function with coefficients defined by 
\begin{equation}\label{Lem2Eq2}
E_{ij}:=\partial_{x_j}\(\left|\nabla v\right|^{p-2}\partial_{x_i}v\)-\frac{1}{n}g\delta_{ij},
\end{equation}
where $\delta_{ij}$ stands for the Kronecker symbol.
\end{lemma}

Now, we prove the following additional result:

\begin{lemma}\label{Lem3}
Let $n\ge2$, $p\in\(1,n\)$, $m,q\in\R$, $u\in W^{1,p}_{\loc}\(\R^n\)\cap L^\infty_{\loc}\(\R^n\)$ be a positive, weak solution of \eqref{IntroEq1}, $v$ and $g$ be the functions defined in \eqref{ProofEq1} and \eqref{ProofEq2}, $E$ be the matrix-valued function defined in \eqref{Lem2Eq2}, and $\varphi$ be a smooth function with compact support in $\R^n$. Then
\begin{multline}\label{Lem3Eq1}
\int_{\R^n}\varphi v^{-q}g^m\(\(\frac{np-n+p}{p}-q\)\left|\nabla v\right|^p+\(\frac{p}{n-p}\)^{p-1}\)\\
+nm\int_{\R^n}\varphi v^{-q}g^{m-1}\left|\nabla v\right|^{p-2}\<E\nabla v,\nabla v\>\\
=-\int_{\R^n}v^{1-q}g^m\left|\nabla v\right|^{p-2}\<\nabla v,\nabla\varphi\>.
\end{multline}
\end{lemma}

\proof[Proof of Lemma~\ref{Lem3}]
By testing \eqref{ProofEq2} against the function $\varphi v^{1-q}g^m$ (which belongs to $C^{0,\alpha'}_{\loc}\(\R^n\)\cap W^{1,2}_{loc}\(\R^n\)$ for some $\alpha'\in\(0,1\)$ according to the above-mentioned regularity properties of the function $v$), we obtain
\begin{multline}\label{Lem3Eq2}
\int_{\R^n}\varphi v^{1-q}g^{m+1}-\(q-1\)\int_{\R^n}\varphi v^{-q}g^m\left|\nabla v\right|^p\\
+m\int_{\R^n}\varphi v^{1-q}g^{m-1}\left|\nabla v\right|^{p-2}\<\nabla g,\nabla v\>\\
+\int_{\R^n}v^{1-q}g^m\left|\nabla v\right|^{p-2}\<\nabla v,\nabla\varphi\>=0.
\end{multline}
The formula \eqref{Lem3Eq1} then follows from \eqref{Lem3Eq2} together with the definition of $g$ and the fact that $\partial_{x_j}g=nv^{-1}E_{ij}\partial_{x_i}v$ for all $j\in\left\{1,\dotsc,n\right\}$ (see Ou~\cite{Ou}*{Lemma~2.1~(i)}).
\endproof

Finally, we state the following results obtained by Ou~\cite{Ou}*{Corollary~2.6 and Lemma~2.7}:

\begin{lemma}\label{Lem4}
Let $n\ge2$, $p\in\(1,n\)$, $u\in W^{1,p}_{\loc}\(\R^n\)\cap L^\infty_{\loc}\(\R^n\)$ be a positive, weak solution of \eqref{IntroEq1}, $v$ and $g$ be the functions defined in \eqref{ProofEq1} and \eqref{ProofEq2} and $E$ be the matrix-valued function defined in \eqref{Lem2Eq2}. Then
\begin{enumerate}
\item[(i)]$\displaystyle\<E^2\nabla v,\nabla v\>\le\Tr\(E^2\)\left|\nabla v\right|^2$
\item[(ii)]For each $n\times n$ matrix-valued function $B$, 
$$\Tr\(BE\)\le \Tr\(E^2\)+C\Tr\(BB^t\)$$
for some constant $C=C\(p\)>0$.
\end{enumerate}
In particular, $\Tr\(E^2\)\ge0$. Moreover, $\Tr\(E^2\)=0$ if and only if $E=0$.
\end{lemma}

We are now in position to prove Theorem~\ref{Th}.

\proof[Proof of Theorem~\ref{Th}]
The beginning of the proof follows ideas from Catino, Monticelli and Roncoroni~\cite{CatMonRon} and Ou~\cite{Ou}. We include it for the sake of completeness. Let $u\in W^{1,p}_{\loc}\(\R^n\)\cap L^\infty_{\loc}\(\R^n\)$ be a positive, weak solution of \eqref{IntroEq1}, $v$ and $g$ be the functions defined in \eqref{ProofEq1} and \eqref{ProofEq2}, and $E$ be the matrix-valued function defined in \eqref{Lem2Eq2}. Let $\eta$ be a smooth, nonnegative cutoff function in $\R^n$ such that $\eta\equiv1$ in $B_1\(0\)$, $\eta\equiv0$ in $\R^n\backslash B_2\(0\)$ and $\left|\nabla\eta\right|\le2$ in $B_2\(0\)\backslash B_1\(0\)$. For each $R>1$, let $\eta_R:\R^n\to\R$ be the function defined as $\eta_R\(x\):=\eta\(x/R\)$ for all $x\in\R^n$, so that $\eta_R\equiv1$ in $B_R\(0\)$, $\eta_R\equiv0$ in $\R^n\backslash B_{2R}\(0\)$ and $\left|\nabla\eta_R\right|\le2/R$ in $B_{2R}\(0\)\backslash B_R\(0\)$. Let $\theta>1$ to be chosen large later on. By using \eqref{Lem2Eq1} with $\varphi=\eta_R^\theta$ and $m=-\frac{p-1}{p}+\varepsilon$ together with Lemma~\ref{Lem4}~(i) and the definition of $g$, we obtain that for small $\varepsilon>0$,
\begin{multline}\label{ProofEq3}
\int_{\R^n}\eta_R^\theta v^{-n}g^{-\frac{2p-1}{p}+\varepsilon}\(n\varepsilon\left|\nabla v\right|^p+\(\frac{p}{n-p}\)^{p-1}\)\Tr\(E^2\)\\
\le-\theta\int_{\R^n}\eta_R^{\theta-1}v^{1-n}g^{-\frac{p-1}{p}+\varepsilon}\left|\nabla v\right|^{p-2}\<E\nabla v,\nabla\eta_R\>.
\end{multline}
Observe that 
\begin{equation}\label{ProofEq4}
n\varepsilon\left|\nabla v\right|^p+\(\frac{p}{n-p}\)^{p-1}\ge\frac{p\varepsilon}{p-1}vg
\end{equation}
provided $\varepsilon$ is chosen small enough. For each $\delta>0$, Lemma~\ref{Lem4}~(ii) with 
$B=-\delta^{-1}\eta_R^{-1}\left|\nabla v\right|^{p-2}\nabla\eta_R\otimes\nabla v$
gives
\begin{align}\label{ProofEq5}
&-\int_{\R^n}\eta_R^{\theta-1}v^{1-n}g^{-\frac{p-1}{p}+\varepsilon}\left|\nabla v\right|^{p-2}\<E\nabla v,\nabla\eta_R\>\nonumber\\
&\quad\le C\delta^{-1}\int_{\R^n}\eta_R^{\theta-2}v^{1-n}g^{-\frac{p-1}{p}+\varepsilon}\left|\nabla v\right|^{2p-2}\left|\nabla \eta_R\right|^2\nonumber\\
&\qquad+\delta\int_{\R^n}\eta_R^\theta v^{1-n}g^{-\frac{p-1}{p}+\varepsilon}\Tr\(E^2\)
\end{align}
for some constant $C=C\(p\)>0$. If $\delta$ is chosen small enough, then it follows from \eqref{ProofEq3}, \eqref{ProofEq4} and \eqref{ProofEq5} that
\begin{align}\label{ProofEq6}
&\int_{\R^n}\eta_R^\theta v^{1-n}g^{-\frac{p-1}{p}+\varepsilon}\Tr\(E^2\)\nonumber\\
&\qquad\le C\int_{\R^n}\eta_R^{\theta-2}v^{1-n}g^{-\frac{p-1}{p}+\varepsilon}\left|\nabla v\right|^{2p-2}\left|\nabla \eta_R\right|^2.
\end{align}
for some constant $C=C\(n,p,\varepsilon,\theta\)>0$. By observing that 
\begin{equation}\label{ProofEq7}
\left|\nabla v\right|\le\(\frac{pvg}{n\(p-1\)}\)^{1/p}
\end{equation}
and since $\left|\nabla \eta_R\right|\le2/R$, we obtain
\begin{equation}\label{ProofEq8}
\int_{\R^n}\eta_R^{\theta-2}v^{1-n}g^{-\frac{p-1}{p}+\varepsilon}\left|\nabla v\right|^{2p-2}\left|\nabla \eta_R\right|^2\le CR^{-2}\int_{\R^n}\eta_R^{\theta-2}v^{-\frac{np-3p+2}{p}}g^{\frac{p-1}{p}+\varepsilon}
\end{equation}
for some constant $C=C\(n,p\)>0$. It follows from \eqref{ProofEq6} and \eqref{ProofEq8} that if 
\begin{equation}\label{ProofEq9}
\int_{\R^n}\eta_R^{\theta-2}v^{-\frac{np-3p+2}{p}}g^{\frac{p-1}{p}+\varepsilon}=\smallo\(R^2\)\quad\text{as }R\to\infty
\end{equation}
and we choose $\theta>2$, then 
\begin{equation}\label{ProofEq10}
\int_{\R^n}v^{1-n}g^{-\frac{p-1}{p}+\varepsilon}\Tr\(E^2\)\le0.
\end{equation}
Since $\Tr\(E^2\)\ge0$ with equality if and only if $E=0$, it then follows from \eqref{ProofEq10} that $E\equiv0$ almost everywhere in $\R^n$, which in turn gives
\begin{equation}\label{ProofEq11}
v\(x\):=c_1+c_2\left|x-x_0\right|^{\frac{p}{p-1}}\quad\forall x\in\R^n
\end{equation}
for some $c_1,c_2\in\R$ (see Catino, Monticelli and Roncoroni~\cite{CatMonRon}*{Section~4.1} or Ciraolo, Figalli and Roncoroni~\cite{CirFigRon}*{Section~3.2}). By putting together \eqref{ProofEq1} and \eqref{ProofEq11} and using \eqref{IntroEq1}, we then obtain that the function $u$ is of the form \eqref{IntroEq2}. Therefore, we are left with showing that \eqref{ProofEq9} holds true. We separate two cases:

\smallskip\noindent
{\bf Case $p>\(n+1\)/3$.} We simplify the arguments used by Ou~\cite{Ou} in this case. By observing that
\begin{align*}
\frac{np-3p+2}{p}+\frac{p-1}{p}+\varepsilon&=\frac{np-2p+1}{p}+\varepsilon<\frac{np-n+p}{p},\allowdisplaybreaks\\
0<\frac{p-1}{p}+\varepsilon&<1
\end{align*}
and
\begin{align*}
&n-\min\(\frac{np-3p+2}{p-1},\frac{np-2p+1}{p}+\varepsilon\)
\\
&\qquad=\max\(\frac{3p-n-2}{p-1},\frac{2p-1}{p}+\varepsilon\)<2
\end{align*}
for small $\varepsilon$, we can apply \eqref{Lem1Eq1}, which gives \eqref{ProofEq9}.

\smallskip\noindent
{\bf Case $p_n<p\le\(n+1\)/3$.} In this case, by observing that for small $\varepsilon$,
\begin{equation}\label{ProofEq12}
\(\frac{p}{n-p}\)^{p-1}\le vg\le\frac{n\(p-1\)}{p\varepsilon}\(\varepsilon\left|\nabla v\right|^p+\(\frac{p}{n-p}\)^{p-1}\),
\end{equation}
we obtain 
\begin{multline}\label{ProofEq13}
\int_{\R^n}\eta_R^{\theta-2}v^{-\frac{np-3p+2}{p}}g^{\frac{p-1}{p}+\varepsilon}\\
\le\(\frac{n-p}{p}\)^{\frac{\(p-1\)\(n-3p+2+p\varepsilon\)}{p}}\int_{\R^n}\eta_R^{\theta-2}v^{-\frac{n\(p-1\)}{p}+\varepsilon}g^{\frac{n-2p+1}{p}+2\varepsilon}
\end{multline}
and
\begin{multline}\label{ProofEq14}
\int_{\R^n}\eta_R^{\theta-2}v^{-\frac{n\(p-1\)}{p}+\varepsilon}g^{\frac{n-2p+1}{p}+2\varepsilon}\le\frac{n\(p-1\)}{p\varepsilon}\\
\times\int_{\R^n}\eta_R^{\theta-2} v^{-\frac{np-n+p}{p}+\varepsilon}g^{\frac{n-3p+1}{p}+2\varepsilon}\(\varepsilon\left|\nabla v\right|^p+\(\frac{p}{n-p}\)^{p-1}\).
\end{multline}
On the other hand, by using \eqref{Lem3Eq1}, we obtain 
\begin{align}\label{ProofEq15}
&\int_{\R^n}\eta_R^{\theta-2} v^{-\frac{np-n+p}{p}+\varepsilon}g^{\frac{n-3p+1}{p}+2\varepsilon}\(\varepsilon\left|\nabla v\right|^p+\(\frac{p}{n-p}\)^{p-1}\)\nonumber\\
&\quad=-\(\theta-2\)\int_{\R^n}\eta_R^{\theta-3}v^{-\frac{n\(p-1\)}{p}+\varepsilon}g^{\frac{n-3p+1}{p}+2\varepsilon}\left|\nabla v\right|^{p-2}\<\nabla v,\nabla\eta_R\>\nonumber\\
&\qquad-n\(\frac{n-3p+1}{p}+2\varepsilon\)\int_{\R^n}\eta_R^{\theta-2} v^{-\frac{np-n+p}{p}+\varepsilon}g^{\frac{n-4p+1}{p}+2\varepsilon}\left|\nabla v\right|^{p-2}\nonumber\\
&\hspace{241pt}\times\<E\nabla v,\nabla v\>.
\end{align}
We begin with estimating the first term in the right-hand side of \eqref{ProofEq15}. For each $\delta>0$ and $q>1$, Young's inequality gives
\begin{align}\label{ProofEq16}
&-\int_{\R^n}\eta_R^{\theta-3}v^{-\frac{n\(p-1\)}{p}+\varepsilon}g^{\frac{n-3p+1}{p}+2\varepsilon}\left|\nabla v\right|^{p-2}\<\nabla v,\nabla\eta_R\>\nonumber\\
&\quad\le\frac{1}{q}\delta^{1-q}\int_{\R^n}\eta_R^{\theta-2-q}v^{-\frac{n\(p-1\)}{p}+\varepsilon}g^{\frac{n-2p+1}{p}-q+2\varepsilon}\left|\nabla v\right|^{q\(p-1\)}\left|\nabla\eta_R\right|^q\nonumber\\
&\qquad+\frac{q-1}{q}\delta\int_{\R^n}\eta_R^{\theta-2}v^{-\frac{n\(p-1\)}{p}+\varepsilon}g^{\frac{n-2p+1}{p}+2\varepsilon}.
\end{align}
By using \eqref{ProofEq7} and since $\left|\nabla \eta_R\right|\le2/R$, we obtain
\begin{multline}\label{ProofEq17}
\int_{\R^n}\eta_R^{\theta-2-q}v^{-\frac{n\(p-1\)}{p}+\varepsilon}g^{\frac{n-2p+1}{p}-q+2\varepsilon}\left|\nabla v\right|^{q\(p-1\)}\left|\nabla\eta_R\right|^q\\
\le CR^{-q}\int_{\R^n}\eta_R^{\theta-2-q}v^{-\frac{\(p-1\)\(n-q\)}{p}+\varepsilon}g^{\frac{n-2p+1-q}{p}+2\varepsilon}
\end{multline}
for some constant $C=C\(n,p,q\)>0$. By observing that $n-2p>0$ when $p\le\(n+1\)/3$, we let
$$q:=n-2p+1+2p\varepsilon,$$
so that 
\begin{align}
q&>1,\nonumber\allowdisplaybreaks\\
\frac{n-2p+1-q}{p}+2\varepsilon&=0\label{ProofEq18}\allowdisplaybreaks\\
\frac{\(p-1\)\(n-q\)}{p}-\varepsilon&=\frac{\(p-1\)\(2p-1\)}{p}-\varepsilon\(2p-1\)\nonumber\\
&\in\(0,\frac{np-n+p}{p}\)\label{ProofEq19}
\end{align}
and
\begin{equation}\label{ProofEq20}
n-q-\frac{\(p-1\)\(n-q\)}{p}+\varepsilon=\frac{2p-1}{p}-\varepsilon<2
\end{equation}
provided $\varepsilon$ is chosen small enough. It follows from \eqref{Lem1Eq1}, \eqref{ProofEq18} and \eqref{ProofEq19} that if $\theta$ is chosen large enough and $\varepsilon$ is chosen small enough, then
\begin{equation}\label{ProofEq21}
R^{-q}\int_{\R^n}\eta_R^{\theta-2-q}v^{-\frac{\(p-1\)\(n-q\)}{p}+\varepsilon}g^{\frac{n-2p+1-q}{p}+2\varepsilon}=\smallo\(R^2\)\,\,\text{ as }R\to\infty.
\end{equation}
By choosing $\delta$ small enough (depending on $n$, $p$, $\varepsilon$ and $\theta$) and putting together \eqref{ProofEq14}, \eqref{ProofEq15}, \eqref{ProofEq16} and \eqref{ProofEq17}, we obtain
\begin{multline}\label{ProofEq22}
\int_{\R^n}\eta_R^{\theta-2}v^{-\frac{n\(p-1\)}{p}+\varepsilon}g^{\frac{n-2p+1}{p}+2\varepsilon}=\smallo\(R^2\)+\bigO\bigg(\int_{\R^n}\eta_R^{\theta-2} v^{-\frac{np-n+p}{p}+\varepsilon}\\
\times
g^{\frac{n-4p+1}{p}+2\varepsilon}\left|\nabla v\right|^{p-2}\left|\<E\nabla v,\nabla v\>\right|\bigg)\quad\text{as }R\to\infty.
\end{multline}
For each $\delta>0$, Lemma~\ref{Lem4}~(ii) with
$$B=\delta^{-1}R^{-1}\eta_R^{-2}v^{\frac{n-2p}{p}}g^{\frac{n-3p}{p}+\varepsilon}\left|\nabla v\right|^{p-2}\nabla v\otimes\nabla v$$ 
gives
\begin{align}\label{ProofEq23}
&\int_{\R^n}\eta_R^{\theta-2} v^{-\frac{np-n+p}{p}+\varepsilon}g^{\frac{n-4p+1}{p}+2\varepsilon}\left|\nabla v\right|^{p-2}\left|\<E\nabla v,\nabla v\>\right|\nonumber\\
&\quad\le C\delta^{-1}R^{-2}\int_{\R^n}\eta_R^{\theta-4} v^{-\frac{np-2n+3p}{p}+2\varepsilon}g^{\frac{2n-7p+1}{p}+3\varepsilon}\left|\nabla v\right|^{2p}\nonumber\\
&\quad\quad+\delta R^2\int_{\R^n}\eta_R^\theta v^{1-n}g^{-\frac{p-1}{p}+\varepsilon}\Tr\(E^2\)
\end{align}
for some constant $C=C\(p\)>0$. By using \eqref{ProofEq6}, \eqref{ProofEq8} and \eqref{ProofEq13}, we obtain
\begin{equation}\label{ProofEq24}
R^2\int_{\R^n}\eta_R^\theta v^{1-n}g^{-\frac{p-1}{p}+\varepsilon}\Tr\(E^2\)\le C
\int_{\R^n}\eta_R^{\theta-2}v^{-\frac{n\(p-1\)}{p}+\varepsilon}g^{\frac{n-2p+1}{p}+2\varepsilon}
\end{equation}
for some constant $C=C\(n,p,\varepsilon,\theta\)>0$. On the other hand, by using \eqref{ProofEq7}, we obtain
\begin{multline}\label{ProofEq25}
\int_{\R^n}\eta_R^{\theta-4} v^{-\frac{np-2n+3p}{p}+2\varepsilon}g^{\frac{2n-7p+1}{p}+3\varepsilon}\left|\nabla v\right|^{2p}\\
\le\(\frac{n\(p-1\)}{p}\)^2\int_{\R^n}\eta_R^{\theta-4} v^{-\frac{np-2n+p}{p}+2\varepsilon}g^{\frac{2n-5p+1}{p}+3\varepsilon}.
\end{multline}
By choosing $\delta$ small enough (depending on $n$, $p$, $\varepsilon$ and $\theta$) and putting together \eqref{ProofEq22}, \eqref{ProofEq23}, \eqref{ProofEq24} and \eqref{ProofEq25}, we obtain
\begin{multline}\label{ProofEq26}
\int_{\R^n}\eta_R^{\theta-2}v^{-\frac{n\(p-1\)}{p}+\varepsilon}g^{\frac{n-2p+1}{p}+2\varepsilon}=\smallo\(R^2\)\\
+\bigO\(R^{-2}\int_{\R^n}\eta_R^{\theta-4} v^{-\frac{np-2n+p}{p}+2\varepsilon}g^{\frac{2n-5p+1}{p}+3\varepsilon}\)\quad\text{as }R\to\infty.
\end{multline}
For each $\delta>0$ and $q>1$, Young's inequality gives
\begin{align}\label{ProofEq27}
&R^{-2}\int_{\R^n}\eta_R^{\theta-4} v^{-\frac{np-2n+p}{p}+2\varepsilon}g^{\frac{2n-5p+1}{p}+3\varepsilon}\nonumber\\
&\quad\le\frac{1}{q}\delta^{1-q}R^{-2q}\int_{\R^n}\eta_R^{\theta-2-2q}v^{-a\(n,p,q,\varepsilon\)}g^{b\(n,p,q,\varepsilon\)}\nonumber\\
&\qquad+\frac{q-1}{q}\delta\int_{\R^n}\eta_R^{\theta-2}v^{-\frac{n\(p-1\)}{p}+\varepsilon}g^{\frac{n-2p+1}{p}+2\varepsilon},
\end{align}
where
$$a\(n,p,q,\varepsilon\):=\frac{np-n-q\(n-p\)}{p}-\varepsilon\(q+1\)$$
and
$$b\(n,p,q,\varepsilon\):=\frac{n-2p+1-q\(3p-n\)}{p}+\varepsilon\(q+2\).$$
If we assume that
\begin{equation}\label{ProofEq28}
0<b\(n,p,q,\varepsilon\)<1
\end{equation}
and
\begin{equation}\label{ProofEq29}
a\(n,p,q,\varepsilon\)+b\(n,p,q,\varepsilon\)<\frac{np-n+p}{p},
\end{equation}
and we choose $\theta$ large enough, then it follows from \eqref{Lem1Eq1} that
\begin{equation}\label{ProofEq30}
\int_{\R^n}\eta_R^{\theta-2-2q}v^{-a\(n,p,q,\varepsilon\)}g^{b\(n,p,q,\varepsilon\)}=CR^{n-\min\(\frac{pa\(n,p,q,\varepsilon\)}{p-1},a\(n,p,q,\varepsilon\)+b\(n,p,q,\varepsilon\)\)}
\end{equation}
for some constant $C=C\(n,p,q,\varepsilon\)>0$. If we assume moreover that
\begin{equation}\label{ProofEq31}
\min\(\frac{pa\(n,p,q,\varepsilon\)}{p-1},a\(n,p,q,\varepsilon\)+b\(n,p,q,\varepsilon\)\)>n-2q-2
\end{equation}
and we choose $\delta$ small enough (depending on $n$, $p$, $\varepsilon$ and $\theta$), then it follows from \eqref{ProofEq26}, \eqref{ProofEq27} and \eqref{ProofEq30} that
\begin{equation}\label{ProofEq32}
\int_{\R^n}\eta_R^{\theta-2}v^{-\frac{n\(p-1\)}{p}+\varepsilon}g^{\frac{n-2p+1}{p}+2\varepsilon}=\smallo\(R^2\)\quad\text{as }R\to\infty.
\end{equation}
Then \eqref{ProofEq9} follows from \eqref{ProofEq13} and \eqref{ProofEq32}. Therefore, it remains to show that for small $\varepsilon$, there exists $q>1$ such that \eqref{ProofEq28}, \eqref{ProofEq29} and \eqref{ProofEq31} simultaneously hold true. When $\varepsilon=0$, we can rewrite \eqref{ProofEq28} as
\begin{equation}\label{ProofEq33}
\frac{n-3p+1}{3p-n}<q<\frac{n-2p+1}{3p-n}
\end{equation}
(observe that $3p-n>0$ since $p>p_n>n/3$). By observing that 
$$a\(n,p,q,0\)+b\(n,p,q,0\)=n-2q-2+\frac{1}{p}>n-2q-2,$$
we can rewrite \eqref{ProofEq29} and \eqref{ProofEq31} with $\varepsilon=0$ as
\begin{equation}\label{ProofEq34}
q>\frac{n-3p+1}{2p}
\end{equation}
and
\begin{equation}\label{ProofEq35}
a\(n,p,q,0\)>\frac{\(p-1\)\(n-2q-2\)}{p},\quad\text{i.e. }q<\frac{2\(p-1\)}{n-3p+2},
\end{equation}
respectively (observe that $n-3p+2\ge1$ since $p\le\(n+1\)/3$). By observing that
$$\frac{n-3p+1}{2p}<\frac{n-3p+1}{3p-n}$$
(recall once again that $n/3<p_n<p\le\(n+1\)/3$), we obtain that \eqref{ProofEq33} implies \eqref{ProofEq34}. On the other hand, it is easy to see that \eqref{ProofEq33} and \eqref{ProofEq35} simultaneously hold true for some $q>1$ if and only if 
\begin{equation}\label{ProofEq36}
\max\(\frac{n-3p+1}{3p-n},1\)<\frac{2\(p-1\)}{n-3p+2}.
\end{equation}
A straightforward computation gives that \eqref{ProofEq36} is equivalent to
$$p>\max\(\frac{n+4}{5},\frac{4n+3-\sqrt{4n^2+12n-15}}{6}\)=p_n.$$
By passing to the limit as $\varepsilon\to0$, we then obtain that if $p>p_n$ and $\varepsilon$ is small enough, then there exists $q>1$ such that \eqref{ProofEq28}, \eqref{ProofEq29} and \eqref{ProofEq31} simultaneously hold true. This ends the proof of Theorem~\ref{Th}.
\endproof

\end{document}